\documentclass[10pt,a4paper,leqno]{amsart}
\usepackage[english]{babel}
\usepackage[T1]{fontenc}
\usepackage{amssymb}
\usepackage{amsfonts}
\usepackage{amsmath}
\usepackage{mathptm}

\swapnumbers
\title{Subharmonic functions, mean value inequality, boundary behavior, nonintegrability and exceptional sets}
\author{Juhani Riihentaus}
\date{\noindent 8 September 2003}
\begin{document}

\maketitle
 
\setcounter{page}{1}

\vspace*{0.1cm}
\noindent\begin{abstract}
\noindent{\footnotesize{
We begin by shortly recalling a generalized mean value inequality for subharmonic functions,  and two applications of it: 
first a  weighted boundary behavior result (with some new references and remarks), and then a borderline
 case result to Suzuki's nonintegrability results for superharmonic and subharmonic
 funtions. The main part of the talk consists, however, of  partial improvements  to Blanchet's removable 
singularity results  for subharmonic, plurisubharmonic and convex functions.}}
\end{abstract}


\setcounter{equation}{0}

\section{Introduction}

\subsection{} In section~2 and 3 we give refinements (Theorems~1 and 2) to our previous results concerning   generalized mean value inequalities for 
subharmonic functions and its applications on the boundary behavior. In section 4 we remark 
that there exists a limiting case result (Theorem~3 and Corollary~3) for Suzuki's results on the 
nonintegrability of  superharmonic
and subharmonic functions. The main part of the article is, however, section 5, where we give partial improvements to 
Blanchet's 
removable singularity results for subharmonic, plurisubharmonic and convex functions (Theorems~4, 5 and Corollaries~4, 5 
and 6). 

\subsection{Notation} Our notation is more or less standard, see
[Ri99]. However, for the convenience of the reader we recall here the following.
We use the common convention $0\cdot \infty =0$. $B(x,r)$ is the Euclidean ball in ${\mathbb{R}}^n$ with center $x$
 and radius $r$. 
We write $\nu_n=m(B(0,1))$, where $m$ is the Lebesgue measure in 
${\mathbb{R}}^n$. In integrals we will  write also $dx$ for the Lebesgue measure. We identify 
${\mathbb{C}}^n$ with ${\mathbb{R}}^{2n}$.
Let $0\leq \alpha \leq n$ and $A\subset {\mathbb{R}}^n$, 
$n\geq 1$. Then we write ${\mathcal{H}}^{\alpha}(A)$ for the
$\alpha$-dimensional Hausdorff (outer) measure of $A$. Recall that $\mathcal{H}^0(A)$ is the number of points of $A$.
In sections 2, 3 and 4  $\Omega$  \emph{is always a domain in} ${\mathbb R}^n$,
$\Omega \ne {\mathbb{ R}}^n$, $n\geq 2$. In section 5 $\Omega$ \emph{is either a domain in} ${\mathbb{R}}^n$ \emph{or in} 
${\mathbb{C}}^n$,
$n\geq 2$. The diameter of $\Omega$ is denoted by 
 diam\,$\Omega$.
The distance from $x\in \Omega$ to $\partial \Omega$, the boundary of $\Omega$,
is denoted by $\delta (x)$. ${\mathcal L}^p_{\mathrm{loc}}(\Omega )$, $p>0$, is the space of
functions $u$ in $\Omega$ for which  $\vert u\vert^p$ is locally integrable on $\Omega$. Our constants $C$ are always
positive, mostly $\geq 1$, and they may vary from line to line. If $x=(x_1,\dots ,x_n)\in {\mathbb R}^n$, $n\geq 2$, and $j\in {\mathbb{N}}$, $1\leq j\leq n$, then we write $x=(x_j,X_j)$, where 
$X_j=(x_1,\dots ,x_{j-1},x_{j+1},\dots ,x_n)$. Moreover, if $A\subset {\mathbb{R}}^n$, $1\leq j\leq n$, and 
$x_j^0\in {\mathbb{R}}$, 
$X_j^0\in {\mathbb{R}}^{n-1}$, we write
\[ A(x_j^0)=\{\, X_j\in {\mathbb{R}}^{n-1}\, :\, \, x=(x_j^0,X_j)\in A\, \},\, \, \, 
 A(X_j^0)=\{\, x_j\in {\mathbb{R}}\, :\, \, x=(x_j,X_j^0)\in A\, \}.\]
We will use similar notation in ${\mathbb{C}}^n$, $n\geq 2$, when considering separately subharmonic and plurisubharmonic 
functions. 

For the definition and properties of subharmonic, separately subharmonic, plurisubharmonic and convex functions, 
see e.g. [Ra37], 
[Le69], [Hel69], [Her71], [H\"o94] and [We94].

\section{The mean value inequality}
\subsection{Previous results}
If $u$ is a nonnegative and subharmonic function on $\Omega$, and $p>0$, then
there is a constant $C=C(n,p)\geq 1$ such that
\begin{equation}
u(x)^p\leq \frac{C}{m(B(x,r))}\int_{B(x,r)}u(y)^p\, dm(y)\end{equation}
for all $B(x,r)\subset \Omega$.
 
See [FeSt72, Lemma~2, p. 172], [Ku74, Theorem~1, p. 529], [Ga81, Lemma~3.7, pp. 121-123],
[AhRu93, (1.5), p. 210]. These authors considered only the case when $u=\vert v\vert$ and 
$v$ is a harmonic function.  However, the proofs in [FeSt72] and [Ga81] apply verbatim also in the general
 case of nonnegative subharmonic functions. This was pointed out in [Ri89, Lemma, p. 69],
[Su90, p. 271], [Su91, p. 113],
 [Ha92, Lemma~1, p. 113], [Pa94, p. 18] and [St98, Lemma 3, p. 305]. In [AhBr88, p. 132] it was pointed
out that a modification of the proof in [FeSt72] gives in fact a slightly more general result, see 2.2 below.
 A possibility for an essentially different proof was pointed out already in [To86, pp. 188-190]. Later other different
 proofs were given in [Pa94, p. 18, and Theorem 1, p. 19] (see also [Pa96, Theorem A, p. 15]), 
[Ri99, Lemma 2.1, p. 233], [Ri00] and [Ri01,Theorem, p. 188].
The results in [Pa94], [Ri99], [Ri00]  and [Ri01] hold in fact 
for more general function classes than just for nonnegative subharmonic functions. See 2.2 and 
Corollary~1  below.  Compare also [DBTr84] and [Do88, p. 485]. 

The inequality $(1)$ has many applications. Among others, it has been
applied to
the (weighted) boundary behavior of nonnegative subharmonic functions
[To86, p. 191], [Ha92,
Theorems~1 and 2, pp. 117-118], [St98, Theorems 1, 2 and 3, pp. 301, 307],
[Ri99, Theorem, p. 233], [Ri00],  and on the nonintegrability of
subharmonic and superharmonic functions
[Su90, Theorem
2, p.
271],
[Su91, Theorem,
 p.
113]. 

Because of the importance of (1), it is worthwhile to present a unified result 
which contains this mean value inequality and all its above referred generalizations.
Such a generalization is proposed below in Theorem~1. In order to state our result and unify the terminology, we give first
two definitions.

\subsection{Quasi-nearly subharmonic functions}
  We call a
(Lebesgue) measurable function
$u:\Omega
\to
 [-\infty ,\infty ]$
 {\textit {quasi-nearly subharmonic}}, if $u\in {\mathcal L}^1_{{\textrm{loc}}}(\Omega
 )$ and if there is a constant $C_0=C_0(n,u,\Omega )\geq 1$ such that
\begin{equation}
u(x)\leq \frac{C_0}{r^n}\int_{B(x,r)}u(y)\, dm(y)\end{equation}
for any ball $B(x,r)\subset \Omega$. Compare [Ri99, p. 233]
and [Do57, p. 430]. Nonnegative quasi-nearly subharmonic
functions have previously been considered in [Pa94]
(Pavlovi\'c called them "functions satisfying the sh$_{\mathrm{K}}$-condition") and in [Ri99], [Ri00] 
(where they were called "pseudosubharmonic functions").
See [Do88, p. 485] for an even  more general function class
of (nonnegative) functions. As a matter of fact, also we will restrict ourselves
 to nonnegative functions.

Nearly
subharmonic functions, thus also quasisubharmonic and subharmonic functions, are examples
of quasi-nearly subharmonic functions.
Recall that a
function $u\in {\mathcal L}^1_{\textrm{loc}}(\Omega )$ is \emph{nearly subharmonic}, if $u$
satisfies $(2)$ with $C_0=\frac{1}{\nu_n}$,
see [Her71, pp. 14, 26]. Furthermore, if $u\geq 0$ is
subharmonic and $p>0$, then by $(1)$ above, $u^p$ is quasi-nearly subharmonic.
By [Pa94, Theorem~1, p. 19] or [Ri99, Lemma 2.1, p. 233] this holds even if $u\geq 0$
is quasi-nearly subharmonic. See also [AhBr88, p. 132].

\subsection{Permissible functions} 

A function \mbox{$\psi :{\mathbb {R}}_+\to {\mathbb {R}}_+$} is called {\textit {permissible}},
if
there  is a nondecreasing, convex function
$\psi_1
:{\mathbb
{R}}_+\to
{\mathbb
{R}}_+$
and an increasing surjection
$\psi_2
:{\mathbb
{R}}_+\to
{\mathbb
{R}}_+$ such that $\psi =\psi_2\circ \psi_1$ and such that the following
conditions are satisfied:
\begin{itemize}
\item[(a)] $\psi_1$ satisfies the
$\varDelta_2$-condition.
\item[(b)] $\psi_2^{-1}$ satisfies the
$\varDelta_2$-condition.
\item[(c)] The function $t\mapsto \frac{t}{\psi_2(t)}$ is
{\textit {quasi-increasing}}, i.e. there is a constant $C=C(\psi_2)\geq 1$ such
that
\begin{displaymath}\frac{s}{\psi_2(s)}\leq C\,    \frac{t}{\psi_2(t)}\end{displaymath}
for all $s, t\in {\mathbb R}_+$, $0\leq s\leq t$.
\end{itemize}
Observe that the condition (b) is equivalent with the following condition.
\begin{itemize}
\item[(b')] For some constant $C=C(\psi_2)\geq 1$,
\begin{displaymath}\psi_2(Ct)\geq 2\, \psi_2(t)\end{displaymath}
for all $t\in {\mathbb R}_+$.
\end{itemize}

Recall that a function $\psi :{\mathbb R}_+\to {\mathbb R}_+$ satisfies the
$\varDelta_2$-{\textit {condition}}, if there is a constant $C=C(\psi )\geq 1$
such that
\begin{displaymath}\psi (2t)\leq C\, \psi (t)\end{displaymath} for all $t\in {\mathbb R}_+$.

If $\psi : {\mathbb
R}_+\to {\mathbb R}_+$ is an 
increasing surjection satisfying the conditions (b) and (c), we say that it is {\textit{strictly permissible}}. Permissible functions are necessarily continuous.

Let it be noted that the condition (c) above is indeed 
natural. For just one counterpart
to it,  see e.g.  [HiPh57, Theorem 7.2.4, p. 239].

Observe that our previous definition for permissible functions in [Ri99, 1.3, p. 232] was {\emph{much more restrictive}}:
 A function $\psi :{\mathbb{R}}_+\rightarrow {\mathbb{R}}_+$ was there 
 defined to be {\emph{permissible}} if it is of the form $\psi (t)=\vartheta (t)^p$,
 $p>0$, where $\vartheta :{\mathbb{R}}_+\rightarrow {\mathbb{R}}_+$ is a nondecreasing, convex
function satisfying the $\varDelta_2$-condition.

\subsection{Examples of permissible functions}  The simple
 example below in (vi), shows that  functions of type (ii)
 are by no means the only permissible functions. The variety  of permissible functions is of course wide.   
\begin{itemize}
\item[(i)]  The functions $\psi_1(t)=\vartheta (t)^p$, $p>0$.
\item[(ii)]  Functions of the form $\psi_2=\phi_2 \circ \varphi_2$,
where
$\phi_2:{\mathbb R}_+\to {\mathbb R}_+$ is a concave surjective function whose
inverse $\phi_2^{-1}$ satisfies the $\varDelta_2$-condition, and
$\varphi_2:{\mathbb R}_+\to {\mathbb R}_+$ is a nondecreasing convex
function
 satisfying the
$\varDelta_2$-condition. (Observe here that any concave function
$\phi_2:{\mathbb R}_+\to {\mathbb R}_+$ is necessarily
nondecreasing.) 
\item[(iii)]
 $\psi_3(t)=c\, t^{p\,\alpha}[\log (\delta +t^{p\,\gamma})]^{\beta}$, where
$c>0$, $0<\alpha <1$, $\delta \geq 1$, and $\beta ,\gamma \in {\mathbb {R}}$
are such that $0<\alpha +\beta \gamma <1$, and $p\geq 1$.
\item[(iv)]  For $0<\alpha <1$, $\beta \geq 0$ and $p\geq 1$,
\begin{displaymath}\psi_4(t)=\left\{ \begin{array}{ll}p^{\beta}t^{p\, \alpha},&{\textrm{ for }}
 0\leq t\leq e,\\
t^{p\, \alpha}(\log t^p)^{\beta},&{\textrm{ for }} t>e.\end{array}\right. \end{displaymath}
\item[(v)] For $0<\alpha <1$, $\beta < 0$ and $p\geq 1$,
\begin{displaymath}\psi_5(t)= \left \{ \begin{array}{ll}
 (\frac{-\beta \,p}{\alpha})^{\beta} t^{p\,\alpha},
&{\textrm{ for }} 0\leq t\leq e^{-{\beta} /{\alpha}},\\
 t^{p\,\alpha}(\log t^p)^{\beta},
&{\text{ for }}
 t>e^{-\beta /{\alpha}}.\end{array} \right. \end{displaymath}
\item[(vi)] For  $p\geq 1$,
\begin{displaymath}\psi_6(t)= \left \{ \begin{array}{lll}
 2n+\sqrt{t^p-2n},&{\textrm{ for }} t^p\in [2n,2n+1), &n=0,1,2,\ldots,\\
 2n+1+[t^p-(2n+1)]^2,
&{\text{ for }}
 t^p\in [2n+1,2n+2), &n=0,1,2,\ldots .\end{array} \right. \end{displaymath}
\end{itemize}

For $p=1$ the functions in (i), (iii), (iv), (v), (vi), and also in (ii) provided that $\varphi_2(t)=t$,
 are 
strictly permissible.

Observe that our previous results were restricted  to the cases 
where $\psi$ was either of \mbox{type (i)} ([Ri99, (1.3), p. 232, and Lemma 2.1, p. 233])  or of
 type (ii) ([Ri01, Theorem, p. 188]).

\subsection{The generalized mean value inequality}  The result (which was presented also at the NORDAN 2000 Meeting, see [Ri00]) 
is the following. Its proof is a 
modification of Pavlovi\'c's argument [Pa94, proof of
Theorem~1, p.
20].
\vskip0.3cm

\noindent\textbf{Theorem~1.} {\emph{Let $u$ be a  nonnegative quasi-nearly subharmonic function on $\Omega$. 
If $\psi :{\mathbb R}_+\to {\mathbb R}_+$ is a permissible function, then $\psi \circ u$ is
quasi-nearly subharmonic on $\Omega$, i.e. there  exists a constant $C=C(n,\psi ,u)\geq 1$
such that
\[ \psi (u(x_0))\leq \frac{C}{{\varrho}^n}\int_{B(x_0,\varrho )}\psi (u(y))\,
dm(y)\] for any ball $B(x_0,\varrho )\subset \Omega$.}}

{\textit {Proof}}. In view of [Ri99, Lemma 2.1, p. 233] we may restrict us to the case where
$\psi =\psi_2:{\mathbb R}_+\to {\mathbb R}_+$ is strictly permissible. 
Since $\psi$ is  continuous, $\psi \circ u$ is measurable and 
$\psi \circ u\in  {\mathcal{L}}^1_{\textrm{loc}}(\Omega)$. It remains to show that $\psi \circ u$ satisfies 
the generalized mean value inequality (2).
But this can be seen exactly as in [Ri01, proof of Theorem, pp. 188-189], the only 
difference being that  instead of the property 2.4 in 
[Ri01, p. 188] of concave functions, one now  uses  the above property (c) in 2.3 
of permissible functions.
\null{} \quad \hfill $\square$

\vskip0.3cm

\noindent\textbf{Corollary~1.} ([Ri01, Theorem, p. 188]){\emph{ Let
$u$ be a  nonnegative subharmonic function on $\Omega$. Let  $\psi :{\mathbb R}_+\to 
{\mathbb R}_+$ be a concave surjection whose inverse $\psi^{-1}$  satisfies the
$\varDelta_2$-condition. Then there exists a constant $C=C(n,\psi ,u)\geq 1$
such that
\[ \psi (u(x_0))\leq \frac{C}{{\varrho}^n}\int_{B(x_0,\varrho )}\psi (u(y))\,
dm(y)\] for any ball $B(x_0,\varrho )\subset \Omega$.}}

\newpage

\section{Weighted boundary behavior}

\subsection{} Before giving our first application of Theorem~1, we recall some terminology from [Ri99, pp.
231--232].

\subsection{Admissible functions} A function
$\varphi :{\mathbb {R}}_+\to {\mathbb {R}}_+$ is {\it {admissible}}, if
it is increasing (strictly), surjective, and  there are constants
$C_2>1$ and
$r_2>0$ such that
\begin{equation*} \varphi
(2t)\leq C_2\, \varphi
(t)
\, \, \, \, 
{\textrm{ {and}}}\, \, \, \, \, \,  \varphi ^{-1} (2s)\leq C_2\, \varphi ^{-1}(s)
\, \, \, \, {\textrm  {for all}}
\, \, \, \, \, \, s, \, t\in {\mathbb {R}}_+,\, \, 0\leq s,\, t\leq
r_2. \end{equation*}

Nonnegative,  nondecreasing functions $\varphi_1(t)$ which satisfy
 the $\varDelta_2$-condition and for which the functions
$t\mapsto \frac{\varphi_1(t)}{t}$ are nondecreasing,
are examples of admissible functions.
Further examples are
$\varphi_2(t)=c\, t^{\alpha}[\log (\delta +t^{\gamma})]^{\beta},$
where
$c>0$,
$\alpha
>0$,
$\delta
\geq
1$,  and
$\beta
,\gamma
\in
{\mathbb
{R}}$  are such that
$\alpha +\beta \gamma >0$.

\subsection{Accessible boundary points  and approach
regions} Let $\varphi :{\mathbb {R}}_+\to {\mathbb {R}}_+$ be an
admissible function and let $\alpha >0$. We say that $\zeta \in \partial
\Omega$ is
$(\varphi ,\alpha )$-{\it {accessible}}, if
\begin{equation} \Gamma_{\varphi}(\zeta ,\alpha )\cap B(\zeta ,\rho )\ne
\emptyset \end{equation} for all
$\rho
 >0$. Here
\[ \Gamma_{\varphi} (\zeta ,\alpha )=\{\, x\in \Omega \, :\, \varphi (\vert
x-\zeta \vert )<\alpha \, \delta (x)\, \},\]
and we call it a  $(\varphi ,\alpha )$-{\textit {approach region
in}}
$\Omega$ {\textit {at}} $\zeta$.
\vskip0.3cm
\subsection{Remarks} (a) In the case when $\varphi (t)=t$, the condition (3) is often called the \emph{corkscrew condition}. 
See e.g. [JeKe82, p. 93].

(b) It follows from [N\"äV\"a91, 2.19, p. 14] that {\emph{all}} boundary points of {\emph{any}} John domain are 
$(\varphi ,\alpha )$-accessible for some $\alpha >0$ (where $\alpha$ depends of course of the parameters of the John domain), 
provided 
the admissible function $\varphi$ satisfies an additional condition,
\begin{equation} \sup \{ \, \frac{\varphi (t)}{t}\, :\, 0<t<r_2\, \}<\infty .\end{equation}
Recall that bounded NTA domains, bounded $(\epsilon , \delta )$-domains of Jones, 
and more generally uniform domains are John domains, see [N\"aV\"a91] and the references therein.
Therefore, using different admissible functions one obtains various types of approach, and in certain cases also non-tangential approach, 
see [St98, pp. 302--304].
 Examples of admissible functions satisfying this additional condition $(4)$ are  
nonnegative,  nondecreasing functions $\varphi_1(t)$ which satisfy
 the $\varDelta_2$-condition and for which the functions
$t\mapsto \frac{\varphi_1(t)}{t}$ are nondecreasing (for small arguments).
Further examples are
$\varphi_2(t)=c\, t^{\alpha}[\log (\delta +t^{\gamma})]^{\beta},$
where
$c>0$,
$\alpha
>1$,
$\delta
\geq
1$,  and
$\beta
,\gamma
\in
{\mathbb
{R}}$  are such that
$\alpha -1 +\beta \gamma >0$.

(c) Mizuta [Mi91] has considered boundary limits of harmonic functions in
Sobolev-Orlicz classes on bounded Lipschitz domains $U$ of ${\mathbb R}^n$, $n\geq
2$. His approach regions are of the form
\[ \Gamma_{\phi} (\zeta ,\alpha )=\{\, x\in U \, :\, \phi (\vert
x-\zeta \vert )<\alpha \, \delta (x)\, \},\]
where now $\phi :{\mathbb R}_+\to {\mathbb R}_+$ is a nondecreasing function which
satisfies the $\varDelta_2$-condition and is such that $t\mapsto \frac{\phi
(t)}{t}$ is nondecreasing. As pointed out above, such functions are admissible
in our sense, and they satisfy also the above condition (4). In fact, they form a proper subclass of our admissible functions.
\vskip0.3cm

\subsection{The weighted boundary behavior result.} Below is the refinement  to our previous
result [Ri99, Theorem, p. 233]. This result was presented also at the NORDAN 2000 Meeting [Ri00], and it improves  the previous 
results of Gehring [Ge57, 
Theorem~1, p. 77], Hallenbeck [Ha92, Theorems~1 and 2, pp. 117-118] and  Stoll [St98,
Theorem~2, p. 307].
\vskip0.3cm

\noindent\textbf{Theorem~2.}
{\emph{Let ${\mathcal{H}}^d(\partial
\Omega )<\infty$ where $0\leq d\leq n$. Suppose that
$u$ is a nonnegative quasi-nearly subharmonic function in $\Omega$.  Let
$\varphi
:{\mathbb {R}}_+\to {\mathbb {R}}_+$ be an admissible function and $\alpha >0$. Let
$\psi
:{\mathbb
{R}}_+\to {\mathbb {R}}_+$ be a permissible function. Suppose that
\begin{equation*} \int_{\Omega}\psi (u(x))\delta (x)^{\gamma}\, dm(x)<\infty \end{equation*}
for some $\gamma \in {\mathbb {R}}$.
Then
\[ \lim_{\rho\to 0}( \sup_{x\in \Gamma_{\varphi ,\rho}(\zeta , \alpha )}\{
\delta
(x)^{n+\gamma}\,
[\varphi
^{-1}(\delta
(x))]^{-d}\,
\psi
(u(x))\}
)
=0\] for ${\mathcal{H}}^d$-almost every $(\varphi ,\alpha )$-accessible point  $\zeta \in
\partial
\Omega$. Here
\[ \Gamma_{\varphi ,\rho}(\zeta ,\alpha )=\{\, x\in \Gamma_{\varphi}(\zeta
,\alpha )\, :\, \delta (x) <\rho \, \}.\]}}

The proof is verbatim the same as [Ri99, proof of
Theorem, pp.
235--238], except that now we just replace [Ri99, Lemma~2.1,
p. 233] 
by the more general Theorem~1 above.
\null{} \quad \hfill \qed
\vskip0.3cm
{\textbf{Remark.}} (Added in December 2003) Mizuta  has given a similar result (for the case when $\psi (t)=t^p$, $p>0$,
 and $\varphi (t)=t^q$, $q\geq 1$) with a different proof,  see 
[Mi01, Theorem 2, p. 73].
\vskip0.3cm 
\noindent\textbf{Corollary~2A.}
{\emph{Let $\Omega$ be a John domain and let ${\mathcal{H}}^d(\partial
\Omega )<\infty$ where $0\leq d\leq n$. Suppose that
$u$ is a nonnegative quasi-nearly subharmonic function in $\Omega$.  Let
$\varphi
:{\mathbb {R}}_+\to {\mathbb {R}}_+$ be an admissible function satisfying the additional condition}} {\textrm{(4)}} {\emph{above. 
Let $\alpha >0$. Let
$\psi
:{\mathbb
{R}}_+\to {\mathbb {R}}_+$ be a permissible function. Suppose that
\begin{equation*} \int_{\Omega}\psi (u(x))\delta (x)^{\gamma}\, dm(x)<\infty \end{equation*}
for some $\gamma \in {\mathbb {R}}$.
Then
\[ \lim_{\rho\to 0}( \sup_{x\in \Gamma_{\varphi ,\rho}(\zeta , \alpha )}\{
\delta
(x)^{n+\gamma}\,
[\varphi
^{-1}(\delta
(x))]^{-d}\,
\psi
(u(x))\}
)
=0\]
 \nopagebreak
\mbox{for ${\mathcal{H}}^d$-almost every   $\zeta \in
\partial
\Omega$}.}}

The proof follows at once from the fact that all boundary points $\zeta \in \partial \Omega$ 
are $(\varphi ,\alpha )$-accessible, as pointed out above in Remark 3.4 (a).
\null{} \quad \hfill \qed
\vskip0.3cm

\noindent\textbf{Corollary~2B.} ([St98, Theorem~2, p. 307]) {\emph{ Let $f$ be a nonnegative subharmonic function on a domain  $G$ in ${\mathbb R}^n$,
$G\ne {\mathbb R}^n$, $n\geq 2$, with ${\mathcal C}^1$ boundary. Let
\begin{equation} \int_{G}f(x)^p\delta (x)^{\gamma}\, dm(x)<\infty \end{equation}
for some $p>0$ and $\gamma >-1-\beta (p)$. Let $0<d\leq n-1$. Then for each
$\tau
\geq
1$ and $\alpha >0$ ($\alpha >1$ when $\tau =1$), there exists a subset
$E_{\tau}$ of
$\partial
G$ with
${\mathcal{H}}^d(E_{\tau})=0$ such that
\begin{equation*} \lim_{\rho\to 0}\{ \sup_{x\in \Gamma_{\tau ,\alpha ,\rho}(\zeta )}[\delta
(x)^{n+\gamma -\frac{d}{\tau}}f(x)^p]\} =0\end{equation*}
for all $\zeta \in \partial G \setminus E_{\tau}$.}}

\vskip0.3cm
Above, for $\zeta \in \partial G$ and
$\varrho >0$,
\[ \Gamma_{\tau ,\alpha ,\varrho}(\zeta )=
     \Gamma_{\tau ,\alpha }(\zeta )\cap G_{\varrho}, \]
where
\[ \Gamma_{\tau ,\alpha}(\zeta )=\{\, x\in G:\, \vert x- \zeta \vert ^{\tau}
<\alpha \, \delta
(x)
\,\},\quad
G_{\varrho}= \{\, x\in G:\, \delta (x)<\varrho \,\}. \]
Moreover,
$\beta (p)=\max \{\, (n-1)(1-p),0\, \}$. 

Stoll makes the assumption $\gamma
>-1-\beta (p)$ in order to exclude the trivial case $f\equiv 0$. As a matter of
fact, it follows from a result of Suzuki [Su90, Theorem~2, p.
271] that $(5)$ together with the condition $\gamma \leq -1-\beta (p)$ implies
indeed that $f\equiv 0$, provided $G$ is a bounded domain with \mbox{$\mathcal{C}^2$ boundary.}  
Unlike Stoll, we have imposed in Theorem~2 no restrictions on the
exponent $\gamma$ in order to exclude the trivial case $u\equiv 0$. Such possibilities are, however, referred in  
Remark 4.5 below.

\section{A Limiting case result  to nonintegrability results of Suzuki}
\subsection{} As another application of Theorem~1, we give in Corollary~3 
below  a supplement, or a limiting case result, to the following result of Suzuki.
\vskip0.3cm

\noindent\textbf{Suzuki's theorem.} ([Su91, Theorem and its proof, pp. 113--115]) 
{\emph{Let $0<p\leq 1$. If a superharmonic (respectively nonnegative subharmonic) function
$v$ on  $\Omega$ satisfies
\begin{equation*} \int_{\Omega}\vert v(x)\vert ^p\, \delta (x)^{np-n-2p}\,
 dm(x)<\infty ,\end{equation*}
then $v$ vanishes identically.
}}             
\vskip0.3cm
Suzuki pointed   out  that his result is sharp in the following sense: If $p$, 
$0<p\leq 1$,
is fixed, then the exponent $\gamma =np-n-2p$ cannot be increased. On the other
hand, clearly $-n<\gamma \leq -2$, when $0<p\leq 1$. Since the
class of permissible
functions include, in addition  to the functions $t^p$, $0<p\leq 1$, also a
large
amount of essentially different functions,
one is tempted to ask whether there exists any limiting case result for Suzuki's result, corresponding to 
the case $p=0$. 
To be more precise, one may pose the following  question: 
\vskip0.3cm

{\emph{Let $\Omega$ and $v$ be as above. Let $\gamma \leq -n$
and let $\psi :{\mathbb{R}}_+\rightarrow {\mathbb{R}}_+$ be permissible. Does the condition
\[ \int_{\Omega}\psi (\vert v(x)\vert )\delta (x)^{\gamma}\, dm(x)<\infty ,\]
\indent imply $v\equiv 0$?}}
\vskip0.3cm

Observe that the least severe form of  above integrability condition occurs when  $\gamma =-n$. 

\subsection{} Before giving an answer in Corollary~3, we state a general
result for arbitrary $\gamma \leq -2$, which is, for $-n< \gamma \leq -2$, however, essentially more or less just
 Suzuki's theorem (see Remarks 4.3 (b) below).
Our formulation has, however, the advantage that, unlike Suzuki's result, it contains  a certain limiting case, Corollary~3, too.
\vskip0.3cm

\noindent\textbf{Theorem~3.}
{\emph{Let $\Omega$ be bounded. Let $v$ be
a superharmonic (respectively nonnegative subharmonic) function  on $\Omega$. Let
$\psi :{\mathbb R}_+\to {\mathbb R}_+$ be a strictly permissible function. Suppose
\begin{equation} \int_{\Omega}\psi (\vert v(x)\vert )\delta (x)^{\gamma}\, dm(x)< \infty ,\end{equation}
 where $\gamma\leq -2$ is such that there is a constant
$C=C(\gamma ,n,\psi ,\Omega )>0$
 for which}}
\begin{equation} s^{n+\gamma}\leq \psi (C\, s^{n-2}) \,{\textit{ for all }}\,
s>\frac{1}{{\textrm{diam}}\,\Omega }.\end{equation} {\emph{Then $v$ vanishes identically.}}
\vskip0.3cm
The
proof is merely a slight modification of Suzuki's  argument, combined
with Theorem~1 above and also some additional estimates. For details, see [Ri03].
\qed 

\subsection{Remarks} Next we consider the assumptions in Theorem~3. 
\begin{itemize}
\item[(a)]  The assumption $\gamma \leq -2$ is unnecessary:
 If
$\gamma \in {\mathbb R}$, then  it follows easily  from
$(7)$ and from the property
(c) in 2.3
of  strictly permissible functions that indeed $\gamma \leq -2$.
\item[(b)] Suppose that $-n<\gamma \leq -2$. If, instead of $(7)$, one supposes that
\begin{equation*} s^{n+\gamma}\leq \psi (C\, s^{n-2}) {\textit{ for all }}
s>0,\end{equation*}
then  clearly
\[ \psi(\vert v(x)\vert )\geq C^{-\frac{n+\gamma}{n-2}}\vert v(x)\vert ^{\frac{n+\gamma}{n-2}}\]
for all $x\in \Omega$. Thus $(6)$ implies that 
\[ \int_{\Omega}\vert v(x)\vert ^{\frac{n+\gamma}{n-2}}\delta (x)^{\gamma}\, dm(x)<\infty ,\]
and hence  $v \equiv 0$ by Suzuki's theorem. Recall that here $0<p=\frac{n+\gamma}{n-2}\leq 1$ and $\gamma =\mbox{np-n-2p}$.
Thus Theorem~3, but now the assumption $(7)$ replaced with the aforesaid assumption, is just a restatement of 
Suzuki's theorem for bounded domains.
\item[(c)] If $\gamma \leq -n$, then the condition $(7)$ clearly holds, since $\psi$ is
strictly permissible. This case gives indeed  the already referred limiting  case for Suzuki's result:
\end{itemize}
\vskip0.3cm 
\noindent\textbf{Corollary~3.}
{\emph{Let $\Omega$ be bounded. Let $v$ be
a superharmonic (respectively nonnegative subharmonic) function  on $\Omega$. Let
$\psi :{\mathbb R}_+\to {\mathbb R}_+$ be  any strictly permissible function and
let $\gamma \leq -n$. If
\begin{equation*} \int_{\Omega}\psi (\vert v(x)\vert )\delta (x)^{\gamma}\, dm(x)<\infty
,\end{equation*} then $v$ vanishes identically.}}

\vskip0.3cm
For the proof observe that the condition $(7)$ is indeed satisfied for $\gamma \leq
-n$, since $\Omega$ is bounded and $\psi$ is increasing.
\null{
}\quad \hfill \qed
\vskip0.3cm

\subsection{Remark} The result of Theorem~3 does not, of course,  hold any
more,
if one replaces strictly permissible functions by permissible functions. For a
counterexample, set, say,  $v(x)=\vert x\vert ^{2-n}$, $\psi (t)=t^p$,
 where $\frac{n-1}{n-2}<p<\frac{n}{n-2}$, $\gamma =np-n-2p$ or just $\gamma
>1$. Then clearly
\[ \int_Bv(x)^p\,\delta (x)^{\gamma}\, dm(x)<\infty \]
but $v\not\equiv 0$.
\vskip0.3cm

\subsection{Remark} Provided $\Omega$ is bounded and $\psi$ is strictly
permissible, one can, with the aid of Theorem~3 and Corollary~3,  exclude some trivial cases $u\equiv 0$ from the result of
Theorem~2  by imposing certain restrictions on the exponent $\gamma$.
 We point out only two cases:
\begin{itemize}
\item[(i)] By Corollary~3, $\gamma >-n$, regardless of $\psi$.
\nopagebreak[4]
\item[(ii)]  By Suzuki's theorem,  $\gamma >np-n-2p$, in
the case when $\psi (t)=t^p$, $0<p\leq 1$.
\end{itemize}

\section{Exceptional sets for subharmonic, plurisubharmonic and convex functions}
 
\subsection{Previous results} Blanchet 
 [Bl95, Theorems 3.1, 3.2 and 3.3, pp. 312--313] gave the following removability results.
\vskip0.3cm

\noindent \textbf{Blanchet's theorem.} {\emph{Let
 $\Omega$ be a domain in ${\mathbb{R}}^n$, $n\geq 2$, and let $S$ be a hypersurface of class ${\mathcal{C}}^1$ which 
divides $\Omega$ into  two subdomains $\Omega_1$ and $\Omega_2$. Let $u\in {\mathcal{C}}^0(\Omega )\cap 
{\mathcal{C}}^2(\Omega_1\cup \Omega_2)$ be subharmonic (respectively convex (or respectively plurisubharmonic 
provided $\Omega$ is then a domain in ${\mathbb{C}}^n$, $n\geq 1$)). If $u_i=u\vert \Omega_i\in {\mathcal{C}}^1
(\Omega_i\cup S)$,
 $i=1,2$, and
\begin{equation} \frac{\partial u_i}{\partial \overline{n}^k}\geq \frac{\partial u_k}{\partial \overline{n}^k}\end{equation}
on $S$ with $i,k=1,2$, then $u$ is subharmonic (respectively convex (or respectively plurisubharmonic)) in $\Omega$.}}
\vskip0.3cm

Above $\overline{n}^k=(\overline{n}_1^k,\dots ,\overline{n}_n^k)$ is the unit normal 
exterior to $\Omega_k$, and $u_k\in {\mathcal{C}}^1(\Omega_k\cup S)$, $k=1,2$, means that there exist $n$ functions $v^j_k$,
 $j=1,\dots ,n$, continuous on $\Omega_k\cup S$, such that 
\[ v_k^j(x)=\frac{\partial u_k}{\partial x_j}(x)\]
for all  $x\in \Omega_k$,  $k=1,2$ and  $j=1, \dots ,n$.

Instead of hypersurfaces of class ${\mathcal{C}}^1$, we will below allow arbitrary sets of finite \mbox{$(n-1)$}-dimensional 
(respectively $(2n-1)$-dimensional) Hausdorff measure as exceptional sets.  Then we must, however, replace the condition 
$(8)$ by another, related condition,  the condition (iv) in Theorem~4 below. In the case of subharmonic and 
plurisubharmonic functions, we  must also impose an additional  integrability
condition on the second partial derivatives {\Large{$\frac{\partial^2 u}{\partial x_j^2}$}}, $j=1,\dots ,n$. 
Observe that in the case of (separately) convex functions we do not, unlike Blanchet,  need any smoothness assumptions of the 
functions (except continuity). Our method of proof is rather elementary, thus natural, with the only exception that we need one geometric 
measure theory result of Federer.

\subsection{The case of subharmonic functions} 
The following  measure theoretic result is essential:
\vskip0.3cm

\noindent \textbf{Lemma~1.} ([Fe69, Theorem 2.10.25, p. 188]) {\emph{Suppose that $A\subset {\mathbb{R}}^n$, $n\geq 2$, is such that
 ${\mathcal{H}}^{n-1}(A)<\infty$. 
Then for all $j$, $1\leq j\leq n$, and for ${\mathcal{H}}^{n-1}$-almost all $X_j\in {\mathbb{R}}^{n-1}$ the set 
$A(X_j)$ is finite.}}
\vskip0.3cm
Our result is:
\vskip0.3cm

\noindent \textbf{Theorem~4.}
{\emph{Suppose that  $\Omega$ is a domain in ${\mathbb{R}}^n$ (respectively in ${\mathbb{C}}^n$), $n\geq 2$. 
Let $E\subset \Omega$ be closed in $\Omega$ 
and ${\mathcal{H}}^{n-1}(E)<\infty$ (respectively ${\mathcal{H}}^{2n-1}(E)<\infty$). Let 
$u:\Omega \to {\mathbb{R}}$ be such that 
\begin{itemize}
\item[(i)] $u\in {\mathcal{C}}^0(\Omega )$,
\item[(ii)] $u\in {\mathcal{C}}^2(\Omega \setminus E)$,
\item[(iii)] for each $j$, $1\leq j\leq n$ (respectively $1\leq j\leq 2n$), 
{\Large{$\frac{\partial ^2u}{\partial x_j^2}$}}$\in {\mathcal{L}}_{{\mathrm{loc}}}^1
(\Omega )$,
\item[(iv)] for each $j$, $1\leq j\leq n$ (respectively $1\leq j\leq 2n$), and for ${\mathcal{H}}^{n-1}$-almost
 all $X_j\in {\mathbb{R}}^{n-1}$ (respectively  for ${\mathcal{H}}^{2n-1}$-almost all $X_j\in {\mathbb{R}}^{2n-1}$)
such that $E(X_j)$ is finite, one has
\[ \liminf_{\epsilon \to 0+0}\frac{\partial u}{\partial x_j}(x_j^0-\epsilon ,X_j)\leq \limsup_{\epsilon \to 0+0}
\frac{\partial u}{\partial x_j}(x_j^0+\epsilon ,X_j)\]
for each $x_j^0\in E(X_j)$,
\item[(v)] $u$ is subharmonic (respectively separately subharmonic) in $\Omega \setminus E$.
\end{itemize}
Then $u$ is subharmonic (respectively separately subharmonic).}}

{\textit{Proof.}} We consider only the subharmonic case. It is sufficient to show that
\begin{equation*}\int u(x)\, \varDelta \varphi (x)\, dx\geq 0\end{equation*}
for all nonnegative testfunctions $\varphi \in {\mathcal{D}}(\Omega )$. Since $u\in {\mathcal{C}}^2(\Omega \setminus E)$
and $u$ is subharmonic in $\Omega \setminus E$, $\varDelta u(x)\geq 0$ for all $x\in \Omega \setminus E$. Therefore 
the claim  follows if we show that 
\begin{equation*}\int u(x)\, \varDelta \varphi (x)\, dx\geq \int  \varDelta u(x)\, \varphi (x)\, dx.\end{equation*}
For this purpose fix $j$, $1\leq j\leq n$, for a while. By Fubini's theorem,
\begin{equation*}\int u(x)\, \frac{\partial ^2\varphi}{\partial x_j^2}(x)\, dx
= \int \bigl[\int u(x_j,X_j)\frac{\partial ^2 \varphi}{\partial x_j^2}(x_j,X_j)\, dx_j\bigr]\, dX_j.\end{equation*}
Using Lemma~1, assumptions (iii), (iv) and Fubini's theorem, we see that for ${\mathcal{H}}^{n-1}$-almost all
 $X_j\in {\mathbb{R}}^{n-1}$,
\begin{equation} \left\{ \begin{aligned}
E(X_j)&{\textrm{ \, is finite}},\\
\frac{\partial^2 u}{\partial x_j^2}(\, \, & \cdot \, \,  ,X_j)\in {\mathcal{L}}_{\mathrm{loc}}^1(\Omega (X_j)),\\
\liminf_{\epsilon \to 0+0}& \frac{\partial u}{\partial x_j}(x_j^0-\epsilon ,X_j)
\leq \limsup_{\epsilon \to 0+0}\frac{\partial u}{\partial x_j}(x_j^0+\epsilon ,X_j){\textrm{ 
\, for all\, }} x_j^0\in E(X_j).\end{aligned}\right.\end{equation}

Let  $K=spt \varphi$.  Choose a domain $\Omega_1$ such that $K\subset \Omega_1\subset \overline{\Omega_1}\subset \Omega$ 
and $\overline{\Omega_1}$ is compact.  Since $E(X_j)$ is finite, there is $M=M(X_j)\in {\mathbb{N}}$ such that
 $E(X_j)=\{\, x_j^1, \ldots ,x_j^M\, \}$ where $x_j^k<x_j^{k+1}$, $k=1, \ldots ,M-1$. Choose for each $k=1, \ldots ,M$ 
 real numbers $a_k, b_k\in (\Omega \setminus E)(X_j)$ such that $a_k<x_j^k<b_k=a_{k+1}<x_j^{k+1}<b_{k+1}$, $k=1,\ldots ,M-1,$ and 
$a_1,b_M\in (\Omega_1\setminus (E\cup K))(X_j)$. Then 
\begin{equation}\int u(x_j,X_j)\, \frac{\partial ^2\varphi}{\partial x_j^2}(x_j,X_j)\, dx_j=
\sum_{k=1}^M\int_{a_k}^{b_k} u(x_j,X_j)\,\frac{\partial^2 \varphi }{\partial x_j^2}(x_j,X_j)\,  dx_j.\end{equation}
Fix $k$, $1\leq k\leq M$, arbitrarily, and write $a=a_k$, $b=b_k$, $x_j^0=x_j^k$. Then
\begin{equation*}\begin{split}\int_a^b u(x_j,X_j) \, \frac{\partial ^2\varphi}{\partial x_j^2}(x_j,X_j)\, dx_j&=
\int_a^{x_j^0} u(x_j,X_j)\, \frac{\partial ^2\varphi}{\partial x_j^2}(x_j,X_j)\, dx_j+
\int_{x_j^0}^b u(x_j,X_j)\, \frac{\partial ^2\varphi}{\partial x_j^2}(x_j,X_j)\, dx_j\\
&=\lim_{\epsilon \to 0+0}\int_a^{x_j^0-\epsilon} u(x_j,X_j)\, \frac{\partial ^2\varphi}{\partial x_j^2}(x_j,X_j)\, dx_j
+\lim_{\epsilon \to 0+0}\int_{x_j^0+\epsilon}^b u(x_j,X_j)\, \frac{\partial ^2\varphi}{\partial x_j^2}(x_j,X_j)\, dx_j\\
&=\lim_{\epsilon \to 0+0}\bigl[{\Big{\vert}}_a^{x_j^0-\epsilon} u(x_j,X_j)\, \frac{\partial \varphi}{\partial x_j}(x_j,X_j)\,
-\int_a^{x_j^0-\epsilon}\frac{\partial u}{\partial x_j}(x_j,X_j)\, \frac{\partial \varphi}{\partial x_j}(x_j,X_j)\, dx_j
\bigr]+\\
&+\lim_{\epsilon \to 0+0}\bigl[{\Big{\vert}}_{x_j^0+\epsilon}^b u(x_j,X_j)\, \frac{\partial \varphi}{\partial x_j}(x_j,X_j)\,
-\int_{x_j^0+\epsilon}^b\frac{\partial u}{\partial x_j}(x_j,X_j)\, \frac{\partial \varphi}{\partial x_j}(x_j,X_j)\, dx_j\bigr]\\
&=\lim_{\epsilon \to 0+0}{\Big{\vert}}_a^{x_j^0-\epsilon} u(x_j,X_j)\, \frac{\partial \varphi}{\partial x_j}(x_j,X_j)
-\lim_{\epsilon \to 0+0}\int_a^{x_j^0-\epsilon}\frac{\partial u}{\partial x_j}(x_j,X_j)\, \frac{\partial \varphi}
{\partial x_j}(x_j,X_j)\, dx_j+\\
&+\lim_{\epsilon \to 0+0}{\Big{\vert}}_{x_j^0+\epsilon}^b u(x_j,X_j)\, \frac{\partial \varphi}{\partial x_j}(x_j,X_j)\,
-\lim_{\epsilon \to 0+0}\int_{x_j^0+\epsilon}^b\frac{\partial u}{\partial x_j}(x_j,X_j)\, 
\frac{\partial \varphi}{\partial x_j}(x_j,X_j)\, dx_j\\
&=u(b,X_j)\frac{\partial \varphi}{\partial x_j}(b,X_j)-u(a,X_j)\frac{\partial \varphi}{\partial x_j}(a,X_j)+\\
&-\lim_{\epsilon \to 0+0}\bigl[{\Big{\vert}}_a^{x_j^0-\epsilon} \frac{\partial u}{\partial x_j}(x_j,X_j)\,\varphi (x_j,X_j)
-\int_a^{x_j^0-\epsilon}\frac{\partial^2 u}{\partial x_j^2}(x_j,X_j)\,  \varphi (x_j,X_j)\, dx_j\bigr]+\\
&-\lim_{\epsilon \to 0+0}\bigl[{\Big{\vert }}_{x_j^0+\epsilon}^{b} \frac{\partial u}{\partial x_j}(x_j,X_j)\,\varphi (x_j,X_j)
-\int_{x_j^0+\epsilon }^{b}\frac{\partial^2 u}{\partial x_j^2}(x_j,X_j)\,  \varphi (x_j,X_j)\, dx_j\bigr].
\end{split}\end{equation*}

Since {\Large{$\frac{\partial^2u}{\partial x_j^2}$}}$(\, \cdot \, , X_j)\in {\mathcal{L}}^1_{\textrm{loc}}(\Omega (X_j))$, the limits
\begin{equation*}\lim_{\epsilon \to 0+0}\int_a^{x_j^0-\epsilon}\frac{\partial^2u}{\partial x_j^2}(x_j,X_j)\, \varphi (x_j,X_j)\, dx_j
{\textrm{ \, and \, }}\lim_{\epsilon \to 0+0}\int_{x_j^0+\epsilon}^b\frac{\partial^2u}{\partial x_j^2}(x_j,X_j)\, 
\varphi (x_j,X_j)\, dx_j\end{equation*}
exist. Thus  also the limits
\begin{equation*}\lim_{\epsilon \to 0+0}{\Big{\vert }}_a^{x_j^0-\epsilon}\frac{\partial u}{\partial x_j}(x_j,X_j)\,
 \varphi (x_j,X_j)
{\textrm{ \, and \, }}\lim_{\epsilon \to 0+0}{\Big{\vert }}_{x_j^0+\epsilon}^b\frac{\partial u}{\partial x_j}(x_j,X_j)\, 
\varphi (x_j,X_j)\end{equation*}
exist. Therefore, remembering also that $a,b\in (\Omega \setminus E)(X_j)$, we get 
\begin{equation}\lim_{\epsilon \to 0+0}{\Big{\vert }}_a^{x_j^0-\epsilon}\frac{\partial u}{\partial x_j}(x_j,X_j)\, 
\varphi (x_j,X_j)=\lim_{\epsilon \to 0+0}\frac{\partial u}{\partial x_j}(x_j^0-\epsilon , X_j)\varphi (x^0_j,X_j)-
\frac{\partial u}{\partial x_j}(a,X_j)\varphi (a,X_j)\end{equation}
and
\begin{equation}\lim_{\epsilon \to 0+0}{\Big{\vert }}_{x_j^0+\epsilon}^b\frac{\partial u}{\partial x_j}(x_j,X_j)\, 
\varphi (x_j,X_j)=\frac{\partial u}{\partial x_j}(b,X_j)\varphi (b,X_j)
-\lim_{\epsilon \to 0+0}\frac{\partial u}{\partial x_j}(x_j^0+\epsilon , X_j)\varphi (x^0_j,X_j).\end{equation}
(The limits 
\begin{equation*}\lim_{\epsilon \to 0+0}\frac{\partial u}{\partial x_j}(x_j^0-\epsilon ,X_j) {\textrm{ \, and \, }}
\lim_{\epsilon \to 0+0}\frac{\partial u}{\partial x_j}(x_j^0+\epsilon ,X_j)\end{equation*}
indeed exist for all points $x_j^0\in (\Omega \setminus E)(X_j)$, for which  $\varphi (x_j^0,X_j)>0$.)

Using $(11)$ and $(12)$ and also the assumption (iv), we get
\begin{equation*}\begin{split}
\int_a^b & u(x_j,X_j)\frac{\partial ^2 \varphi}{\partial x_j^2}(x_j,X_j)\, dx_j
=\bigl[u(b,X_j)\,\varphi (b,X_j)-u(a,X_j)\, \varphi (a,X_j)\bigr]+\\
&+\bigl[-\frac{\partial u}{\partial x_j}(b,X_j)\, \varphi (b,X_j)+\frac{\partial u}{\partial x_j}(a,X_j)\, \varphi (a,X_j)\bigr]
+\bigl[\lim_{\epsilon \to 0+0}\frac{\partial u}{\partial x_j}(x_j^0+\epsilon , X_j)-
\lim_{\epsilon \to 0+0}\frac{\partial u}{\partial x_j}(x_j^0-\epsilon ,X_j)\bigr]\,\varphi (x^0_j,X_j)+\\
&+\lim_{\epsilon \to 0+0}\int_a^{x_j^0-\epsilon }\frac{\partial ^2u}{\partial x_j^2}(x_j,X_j)\, \varphi (x_j,X_j)\, dx_j
+\lim_{\epsilon \to 0+0}\int_{x_j^0+\epsilon}^b\frac{\partial ^2u}{\partial x_j^2}(x_j,X_j)\, \varphi (x_j,X_j)\, dx_j\geq \\
&\geq \bigl[u(b,X_j)\, \varphi (b,X_j)-u(a,X_j)\, \varphi (a,X_j)\bigr]+\bigl[-\frac{\partial u}{\partial x_j}(b,X_j)\,
 \varphi (b,X_j)+\frac{\partial u}{\partial x_j}(a,X_j)\,
 \varphi (a,X_j)\bigr]+\int_a^b \frac{\partial ^2u}{\partial x_j^2}(x_j,X_j)\, \varphi (x_j,X_j)\, dx_j. 
\end{split}\end{equation*}
In view of this and of $(10)$ we get
\begin{equation*} \int u(x_j,X_j)\, \frac{\partial^2\varphi}{\partial x_j^2}(x_j,X_j)\, dx_j\geq
\int  \frac{\partial^2u}{\partial x_j^2}(x_j,X_j)\,\varphi (x_j,X_j)\, dx_j.\end{equation*}
Integrating then here on both sides with respect to $X_j\in {\mathbb{R}}^{N-1},$ and using (9) and also Fubini's theorem, we get
\begin{equation*} \int u(x)\, \frac{\partial ^2\varphi}{\partial x_j^2}(x)\, dx\geq
\int  \frac{\partial^2u}{\partial x_j^2}(x)\,\varphi (x)\, dx.\end{equation*}
Hence
\begin{equation*} \int u(x)\, \varDelta \varphi (x)\, dx_j\geq
\int  \varDelta u(x)\,\varphi (x)\, dx\geq 0,\end{equation*}
concluding the proof. \qed
\vskip0.3cm

\noindent \textbf{Corollary~4A.}
{\emph{Suppose that  $\Omega$ is a domain in ${\mathbb{R}}^n$ (respectively in ${\mathbb{C}}^n$), $n\geq 2$. Let $E\subset \Omega$ be closed in $\Omega$ 
and ${\mathcal{H}}^{n-1}(E)<\infty$ (respectively ${\mathcal{H}}^{2n-1}(E)<\infty$). Let 
$u:\Omega \to {\mathbb{R}}$ be such that 
\begin{itemize}
\item[(i)] $u\in {\mathcal{C}}^1(\Omega )$,
\item[(ii)] $u\in {\mathcal{C}}^2(\Omega \setminus E)$,
\item[(iii)] for each $j$, $1\leq j\leq n$ (respectively $1\leq j\leq 2n$), 
{\Large{$\frac{\partial ^2u}{\partial x_j^2}$}}$\in {\mathcal{L}}_{{\mathrm{loc}}}^1
(\Omega )$,
\item[(iv)] $u$ is subharmonic (respectively separately subharmonic) in $\Omega \setminus E$.
\end{itemize}
Then $u$ is subharmonic (respectively separately subharmonic).}}

\subsection{The case of plurisubharmonic functions}
In order to obtain a similar result for plurisubharmonic functions, we need  the following result of Lelong.

\noindent \textbf{Lemma~2.} ([Le69, Theorem 1, p. 18])
{\emph{ Suppose that $D$ is a domain of ${\mathbb{C}}^n$, $n\geq 2$.
 Let $v:D \to [-\infty ,+\infty )$. Then $v$ is plurisubharmonic if and only if the  following condition holds:\newline 
For each $z_0\in D$ and for each affine transformation 
$A=(A_1, \ldots ,A_n): {\mathbb{C}}^n\to {\mathbb{C}}^n$,
\begin{equation*}\begin{split}
z'&=Az \Leftrightarrow (z'_1,\dots ,z'_n)=(A_1(z_1,\dots ,z_n),\dots ,A_n(z_1,\dots ,z_n))\\
\Leftrightarrow & \left\{ \begin{split}z_1'&=A_1(z_1,\dots ,z_n)=z_1^0+a_{11}z_1+\cdots +a_{1n}z_n,\\
\vdots &\\ 
z_n'&=A_n(z_1,\ldots ,z_n)=z_n^0+a_{n1}z_1+\cdots +a_{nn}z_n,
\end{split}
\right. \end{split}\end{equation*}
for which $\det A\ne 0$, the function $v\circ A: A^{-1}(D)\to [-\infty ,+\infty )$ is subharmonic.}}
\vskip0.3cm

\noindent \textbf{Corollary~4B.}
{\emph{ Suppose that $\Omega$ is a domain of  ${\mathbb{C}}^n$, $n\geq 2$.  
Let $E\subset \Omega$ be closed in $\Omega$ 
and ${\mathcal{H}}^{2n-1}(E)<\infty$. Let 
$u:\Omega \to {\mathbb{R}}$ be such that 
\begin{itemize}
\item[(i)] $u\in {\mathcal{C}}^1(\Omega )$,
\item[(ii)] $u\in {\mathcal{C}}^2(\Omega \setminus E)$,
\item[(iii)] for each $j$, $1\leq j\leq 2n$, 
{\Large{$\frac{\partial ^2u}{\partial x_j^2}$}}$\in {\mathcal{L}}_{{\mathrm{loc}}}^1
(\Omega )$,
\item[(iv)] $u$ is plurisubharmonic in $\Omega \setminus E$.
\end{itemize}
Then $u$ is plurisubharmonic.}}

{\textit{Proof.}} By Lemma~2  it is sufficient to show that $v=u\circ A$ is subharmonic in $\Omega ' =A^{-1}(\Omega )$ 
for
any affine mapping $A: \, {\mathbb{C}}^n\to {\mathbb{C}}^n$ with $\det A\ne 0$.
Clearly $v\in {\mathcal{C}}^1(\Omega ')$ and $v\in {\mathcal{C}}^2(\Omega ' \setminus E')$, where $E'=A^{-1}(E)$.
It is easy to see that for each $j$, $1\leq j\leq 2n$, 
{\Large{$\frac{\partial ^2v}{\partial x_j^2}$}}$\in {\mathcal{L}}^1_{\mathrm{loc}}(\Omega ')$. Since $u$ is 
plurisubharmonic 
in $\Omega \setminus E$, $v$ is by Lemma~2 subharmonic in $\Omega '\setminus E'$, thus subharmonic in $\Omega'$ by 
Corollary~4A.\qed

\subsection{The case of convex functions} We recall first some very basic properties of convex functions.
 
Let $D$ be a domain of of ${\mathbb{R}}^n$, $n\geq 1$. A function $f:D\to {\mathbb{R}}$ is 
\emph{convex} if the following 
condition is satisfied: For each $x,y\in D$ such that $\{ \,tx+(1-t)y\, :\, t\in [0,1]\, \}\subset  D$, one has $f(tx+(1-t)y)\leq t\, f(x)+(1-t)\, f(y)$
for all $t\in [0,1]$. 
\vskip0.3cm

\noindent \textbf{Lemma~3.} ([We94, Theorem~5.1.3, p.~195]) 
{\emph{ Let $f:[a,b]\to {\mathbb{R}}$ be a convex function. Then $f$ possesses left and right derivatives at each interior point 
of $[a,b]$, and if $x_1,x_2$ are interior points of $[a,b]$ with $x_1<x_2$, then   
\[ -\infty < f'_{-}(x_1)\leq f'_{+}(x_1)\leq \frac{f(x_2)-f(x_1)}{x_2-x_1}\leq f'_{-}(x_2)\leq f'_{+}(x_2)<+\infty .\]
}}
\vskip0.3cm

\noindent \textbf{Lemma~4.} ([We94, Theorem 5.1.8, p. 198]) 
{\emph{ Let $f:(a,b)\to {\mathbb{R}}$. Then $f$ is convex if and only if it has support at each point of $(a,b)$, 
i.e. for any $x_0\in (a,b)$ there is a constant $m\in {\mathbb{R}}$ such that
\[f(x_0)+m\, (x-x_0)\leq f(x)\]
for all $x\in (a,b)$.}}
\vskip0.3cm

Moreover, if $f$ is convex, then any  $m$, $f'_-(x_0)\leq m\leq f'_+(x_0)$, will do.
\vskip0.3cm
We consider first separately convex functions:
\vskip0.3cm
 
\noindent\textbf{Theorem~5.}
{\emph{ Suppose that $\Omega$ is a domain of ${\mathbb{R}}^n$, $n\geq 2$.
 Let $E\subset \Omega$ be closed in $\Omega$ and ${\mathcal{H}}^{n-1}(E)<\infty$. Let $u:\Omega \to {\mathbb{R}}$ be 
such that 
\begin{itemize}
\item[(i)] $u\in {\mathcal{C}}^0(\Omega )$,
\item[(ii)] for each $j$, $1\leq j\leq n$, and for ${\mathcal{H}}^{n-1}$-almost all $X_j\in {\mathbb{R}}^{n-1}$
such that $E(X_j)$ is finite, one has
\[\liminf_{\epsilon \to 0+0}\frac{\partial_- u}{\partial x_j}(x_j^0-\epsilon ,X_j)\leq \limsup_{\epsilon \to 0+0}
\frac{\partial_+ u}{\partial x_j}(x_j^0+\epsilon ,X_j)\]
for each $x_j^0\in E(X_j)$,
\item[(iii)] $u$ is separately convex in $\Omega\setminus E$.
\end{itemize}Then $u$ is separately convex.}}

\vskip0.3cm
Above,  and in the sequel, {\Large{$\frac{\partial_-u}{\partial x_j}$}}$(x_j,X_j)$ and 
{\Large{$\frac{\partial_+u}{\partial x_j}$}}$(x_j,X_j)$, $j=1,\ldots ,n$, are the left and right 
partial derivatives of $u$, respectively, taken at the point $x=(x_j,X_j)$.

Observe that the condition (ii) is  \emph{a necessary condition} for (separately) convex 
functions. 
 
{\textit{Proof of Theorem~5.}} Choose $j$, $1\leq j\leq n$, arbitrarily. Using Lemma~1 and the condition (ii) we see 
that for ${\mathcal{H}}^{n-1}$-almost all $X_j\in {\mathbb{R}}^{n-1}$,
\begin{equation*} \left\{ \begin{aligned}
E(X_j)&{\textrm{ \, is finite}},\\
(\Omega \setminus E)&(X_j)\ni x_j\mapsto u(x_j,X_j)\in {\mathbb{R}} {\textrm{ \, is convex}},\\
\liminf_{\epsilon \to 0+0}& \frac{\partial_{-} u}{\partial x_j}(x_j^0-\epsilon ,X_j)
\leq \limsup_{\epsilon \to 0+0}\frac{\partial_{+} u}{\partial x_j}(x_j^0+\epsilon ,X_j) {\textrm{\, for all\, }} x_j^0
\in E(X_j).\end{aligned}\right.\end{equation*}
Using this, Lemma~3 and Lemma~4 one sees that for ${\mathcal{H}}^{n-1}$-almost all  $X_j\in {\mathbb{R}}^{n-1}$ the 
functions
\begin{equation}\Omega (X_j)\ni x_j\mapsto u(x_j,X_j)\in{\mathbb{R}}\end{equation}
are in fact convex. (Here one proceeds  e.g. as follows: Suppose that  $(a,b)$ is an arbitrary interval of the open set 
$\Omega(X_j)$, that $E(X_j)\cap (a,b)=\{\, x_j^1,\dots ,x_j^N\, \}$, where $a<x_j^k<x_j^{k+1}<b$, $k=1, \dots, N-1$ and 
$x_j^{N+1}=b$. If 
$u(\cdot ,X_j)\vert (a,x_j^k)$ and $u(\cdot ,X_j)\vert (x_j^k,x_j^{k+1})$,  are convex, then 
it follows from the assumptions 
that  $u(\cdot ,X_j)\vert (a,x_j^{k+1})$, is convex,  $k=1,\dots ,N$.) From this and from the fact that $u$ is continuous,
 it follows easily that the functions of the  form (13) above  
are in fact convex  
{\emph{for all}} $X_j\in {\mathbb{R}}^{n-1}$. Since $j$, $1\leq j\leq n$, was arbitrary, the claim follows.\qed
\vskip0.3cm

\noindent\textbf{Corollary~5.}
{\emph{ Suppose that $\Omega$ is a domain of ${\mathbb{R}}^n$, $n\geq 2$.
 Let $E\subset \Omega$ be closed in $\Omega$ and ${\mathcal{H}}^{n-1}(E)<\infty$. Let $u:\Omega \to {\mathbb{R}}$ be 
such that 
\begin{itemize}
\item[(i)] $u\in {\mathcal{C}}^1(\Omega )$,
\item[(ii)] $u$ is (separately) convex in $\Omega\setminus E$.
\end{itemize}
Then $u$ is (separately) convex.}}
\vskip0.3cm

The separately convex case follows directly from Theorem~5. The convex case follows from the separately convex case 
with the aid of the following  Lelong type result (whose proof is 
similar to  [Le69, proof of Theorem 1, p. 18]). 
\qed
\vskip0.3cm
\noindent \textbf{Lemma~5.} 
{\emph{ Suppose that $D$ is a domain of ${\mathbb{R}}^n$, $n\geq 2$.
 Let $v:D \to [-\infty ,+\infty )$. Then $v$ is convex if and only if the  following condition holds:\newline 
For each $x_0\in D$ and for each affine transformation 
$A=(A_1, \ldots ,A_n): {\mathbb{R}}^n\to {\mathbb{R}}^n$,
\begin{equation*}\begin{split}
x'&=Ax  \Leftrightarrow (x'_1,\dots ,x'_n)=(A_1(x_1,\dots ,x_n),\dots ,A_n(x_1,\dots ,x_n))\\
\Leftrightarrow & \left\{ \begin{split}x_1'&=A_1(x_1,\dots ,x_n)=x_1^0+a_{11}x_1+\cdots +a_{1n}x_n,\\
\vdots &\\ 
x_n'&=A_n(x_1,\ldots ,x_n)=x_n^0+a_{n1}x_1+\cdots +a_{nn}x_n,
\end{split}
\right. \end{split}\end{equation*}
for which $\det A\ne 0$, the function $v\circ A: A^{-1}(D)\to [-\infty ,+\infty )$ is separately convex.}}

\vskip0.3cm
\bibliographystyle{alpha}
\begin{center}\textsc{References}\end{center}
\footnotesize{
\begin{enumerate}
\item[{[AhBr88]}] P. Ahern and J. Bruna, \emph{Maximal and area integral characterizations of Hardy--Sobolev
spaces in the unit ball of} ${\mathbb{C}}^n$, Revista Matemática Iberoamericana \textbf{4}(1988), 123--153.
\item[{[AhRu93]}] P. Ahern and W. Rudin, \emph{Zero sets of functions in harmonic 
Hardy spaces}, Math. Scand.  \textbf{73}(1993), 209--214.
\item[{[Bl95]}] P. Blanchet, \emph{On removable singularities of subharmonic and plurisubharmonic functions},
Complex Variables  \textbf{26}(1995), 311--322.
\item[{[BlGa94]}] P. Blanchet  and P.M. Gauthier, \emph{Fusion of two solutions of a partial differential equation}, 
Meth, Appl. Anal., \textbf{1}, no. 3(1994), 371--384.
\item[{[DBTr84]}] E. Di Benedetto and N.S. Trudinger, \emph{Harnack inequalities for quasi-minima of variational 
integrals}, Ann. Inst. H. Poincaré, Analyse Nonlineaire  \textbf{1}(1984), 295--308.
\item[{[Do57]}] Y. Domar,\emph{ On the existence of a largest subharmonic minorant of a given
function}, Arkiv f\"or Mat. \textbf{3}(1957), 429--440.
\item[{[Do88]}] Y. Domar, \emph{Uniform boundedness in families related to subharmonic
 functions},  J. London Math. Soc. (2) \textbf{38}(1988), 485--491.
\item[{[Fe69]}] H. Federer, \emph{Geometric Measure Theory},
Springer-Verlag, Berlin, 1969.
\item[{[FeSt72]}] C. Fefferman and E.M. Stein, \emph{H$^p$ spaces of several variables},
Acta Math.  \textbf{l29}(1972), 137--193.
\item[{[Ga81]}] J.B. Garnett, \emph{Bounded analytic functions}, Academic Press, 
 New York, 1981.
\item[{[Ge57]}]  F.W. Gehring, \emph{On the radial order of subharmonic functions}, 
J.\ Math. Soc. Japan \textbf{9}(1957), 77--79.
\item[{[Ha92]}] D.J. Hallenbeck, \emph{Radial growth of subharmonic functions}, 
Pitman Research Notes 262, 1992, pp. 113--121.
\item[{[Hel69]}] L.L. Helms, \emph{Introduction to potential theory}, Wiley-Interscience,
New York, 1969.
\item[{[Her71]}] M. Hervé, \emph{Analytic and Plurisubharmonic Functions in Finite and Infinite
Dimensional Spaces},  Lecture Notes in Mathematics 198,  Springer-Verlag, Berlin, 1971.
\item[{[HiPh57]}] E.~Hille and R.S.~Phillips, \emph{Functional Analysis and Semigroups}, 
American Mathematical Society, Colloquium publications XXXI, Providence, R.I., 1957.
\item[{[H\"o94]}] L.~H\"ormander, \emph{Notions of Convexity}, Birkh\"auser, Boston, 1994.
\item[{[JeKe82]}] D.~Jerison and C.~Kenig, \emph{Boundary behavior of harmonic functions in nontangentially
accessible domains}, Adv. Math. \textbf{46}(1982), 80--147.
\item[{[Kr83]}] J. Král, \emph{Some extension results concerning subharmonic functions}, J. London Math. Soc. 2,
 \textbf{28}(1983), 62--70.
\item[{[Ku74]}] \"U. Kuran, \emph{Subharmonic behavior of 
$\vert h\vert ^p$ ($p>0$, $h$ harmonic)},  J. London Math. Soc. (2)  \textbf{8}(1974),
529--538.
\item[{[Le69]}] P. Lelong, \emph{Plurisubharmonic functions and positive differential forms}, 
Gordon and Breach, New York, 1969.
\item[{[Mi91]}] Y. Mizuta, \emph{Boundary limits of harmonic functions in Sobolev-Orlicz
classes},  Potential Theory (ed. M. Kishi), Walter de Gruyter\&Co, 
 Berlin $\cdot$ New York, 1991, pp. 235--249.
\item[{[Mi01]}] Y.~Mizuta, \emph{Boundary limits of functions in weighted Lebesgue or Sobolev classes},
Rev. Roumaine Math. Pures Appl. \textbf{46}(2001), 67--75.
\item[{[N\"aV\"a91]}] R.~N\"akki and J.~V\"ais\"al\"a, \emph{John disks},  Exposition Math. \textbf{9}(1991), 3--43.
\item[{[Pa94]}] M. Pavlovi\'c, \emph{On subharmonic behavior and oscillation of functions
 on balls in} ${\mathbb R}^n$, 
 Publ. de l'Institut Mathém., Nouv. sér.  \textbf{55(69)}(1994), 18--22.
\item[{[Pa96]}] M. Pavlovi\'c, \emph{Subharmonic behavior of smooth functions},
 Math. Vesnik  \textbf{48}(1996), 15--21.
\item[{[Ra37]}] T. Rado, \emph{Subharmonic Functions}, Springer-Verlag, Berlin, 1937.
\item[{[Ri89]}] J. Riihentaus, \emph{On a theorem of Avanissian--Arsove},
Expo. Math. \textbf{7}(1989), 69--72.
\item[{[Ri99]}] J. Riihentaus, \emph{Subharmonic functions: non-tangential and 
tangential boundary
behavior},  Function Spaces, Differential Operators and Nonlinear Analysis,
Proceedings of the Sy\"ote Conference 1999 (Pudasj\"arvi, Finland, June 10--16, 1999)  (eds.  V. Mustonen and J. R\'akosnik), 
 Math. Inst., Czech Acad. Science,  Praha, 2000, pp. 229--238.
\item[{[Ri00]}] J. Riihentaus, \emph{Subharmonic functions: non-tangential and 
tangential boundary
behavior} (Abstract),  Proceedings of The Nordic Complex  Analysis Meeting NORDAN 2000 (\"Ornsk\"oldsvik, Sweden, May 5--7, 
2000), to appear. 
\item[{[Ri01]}] J. Riihentaus, \emph{A generalized mean value inequality for 
subharmonic functions}, Expo. Math.  \textbf{19}(2001), 187--190.
\item[{[Ri03]}] J. Riihentaus, \emph{A generalized mean value inequality for 
subharmonic functions and applications}, arXiv:math.CA/0302261 v1 21 Feb 2003.
\item[{[St98]}] M. Stoll, \emph{Weighted tangential boundary limits of subharmonic functions on
domains in ${\mathbb {R}}^n$ ($n\geq 2$)}, Math. Scand.  \textbf{83}(1998), 300--308.
\item[{[Su90]}]  N. Suzuki, \emph{Nonintegrability of harmonic functions in a domain},
Japan. J. Math. \textbf{16}(1990), 269--278.
\item[{[Su91]}] N. Suzuki, \emph{Nonintegrability of superharmonic functions}, 
Proc.\ Amer.\ Math.\ Soc. \ \textbf{113}(1991), 113--115.
\item[{[To86]}]  A. Torchinsky, \emph{Real-Variable Methods in Harmonic Analysis},
 Academic Press,  London, 1986.
\item[{[We94]}] R. Webster, \emph{Convexity}, Oxford University Press, Oxford, 1994.
\end{enumerate}
\vskip0.2cm
\footnotesize{
\noindent\textsc{South Carelia Polytechnic}\\
\textsc{P.O. Box 99}\\
\textsc{FIN-53101 Lappeenranta, Finland}\\

\noindent\textsc{and}\\

\noindent\textsc{Department of Mathematics}\\
\textsc{University of Joensuu}\\
\textsc{P.O. Box 111}\\
\textsc{FIN-80101 Joensuu, Finland}\\

\noindent E-mail address: \texttt{juhani.riihentaus{@}scp.fi}
}
\end{document}